\documentclass[11pt,letterpaper]{amsart}

\usepackage[english]{babel}
\usepackage{amsmath}
\usepackage{graphicx}
\usepackage[dvipsnames]{xcolor}
\usepackage[colorlinks, allcolors=MidnightBlue]{hyperref}
\usepackage{pgfplots}
\pgfplotsset{compat=1.15}
\usepackage{mathrsfs}
\usetikzlibrary{arrows}

\usepackage{amssymb,amsmath,amsthm,amscd,mathtools,mathdots}
\usepackage{multirow}
\usepackage{enumitem}
\usepackage{caption}
\usepackage{subcaption}
\usepackage{tikz}
\usepackage{tikz-cd}

\newtheorem{thm}{Theorem}[section]
\newtheorem{lm}[thm]{Lemma}
\newtheorem{coro}[thm]{Corollary}
\newtheorem{prop}[thm]{Proposition}

\theoremstyle{definition}
\newtheorem{df}[thm]{Definition}

\theoremstyle{plain}

\newcommand{\R}{\mathbb{R}}
\newcommand{\C}{\mathbb{C}}
\newcommand{\Z}{\mathbb{Z}}
\newcommand{\A}{\mathbb{A}}
\newcommand{\Au}{\mathbb{A}^1}

\newcommand{\GW}{\operatorname{GW}}

\newcommand{\Pk}{\mathbb{P}_{k}^2}

\newcommand{\val}{\operatorname{val}}

\newcommand{\mult}[1][k]{\operatorname{mult}^{\A^1}_{#1}}

\newcommand{\qinv}[1]{\left\langle#1\right\rangle}

\newcommand{\Tra}[2][E/k]{\operatorname{Tr}_{#1}\left(#2\right)}

\newcommand{\Puiseux}{\{\!\{t\}\!\}}

\newcommand{\sig}[1]{\operatorname{sig}\left(#1\right)}

\newcommand{\Wel}{\operatorname{Wel}}

\newcommand{\Bl}{\operatorname{Bl}_p\left(\Sigma_1\right)}

\makeatletter
\let\save@mathaccent\mathaccent
\newcommand*\if@single[3]{%
    \setbox0\hbox{${\mathaccent"0362{#1}}^H$}%
    \setbox2\hbox{${\mathaccent"0362{\kern0pt#1}}^H$}%
    \ifdim\ht0=\ht2 #3\else #2\fi
}
\newcommand*\rel@kern[1]{\kern#1\dimexpr\macc@kerna}
\newcommand*\widebar[1]{\@ifnextchar^{{\wide@bar{#1}{0}}}{\wide@bar{#1}{1}}}
\newcommand*\wide@bar[2]{\if@single{#1}{\wide@bar@{#1}{#2}{1}}{\wide@bar@{#1}{#2}{2}}}
\newcommand*\wide@bar@[3]{%
    \begingroup
    \def\mathaccent##1##2{%
        \let\mathaccent\save@mathaccent
        \if#32 \let\macc@nucleus\first@char \fi
        \setbox\z@\hbox{$\macc@style{\macc@nucleus}_{}$}%
        \setbox\tw@\hbox{$\macc@style{\macc@nucleus}{}_{}$}%
        \dimen@\wd\tw@
        \advance\dimen@-\wd\z@
        \divide\dimen@ 3
        \@tempdima\wd\tw@
        \advance\@tempdima-\scriptspace
        \divide\@tempdima 10
        \advance\dimen@-\@tempdima
        \ifdim\dimen@>\z@ \dimen@0pt\fi
        \rel@kern{0.6}\kern-\dimen@
        \if#31
        \overline{\rel@kern{-0.6}\kern\dimen@\macc@nucleus\rel@kern{0.4}\kern\dimen@}%
        \advance\dimen@0.4\dimexpr\macc@kerna
        \let\final@kern#2%
        \ifdim\dimen@<\z@ \let\final@kern1\fi
        \if\final@kern1 \kern-\dimen@\fi
        \else
        \overline{\rel@kern{-0.6}\kern\dimen@#1}%
        \fi
    }%
    \macc@depth\@ne
    \let\math@bgroup\@empty \let\math@egroup\macc@set@skewchar
    \mathsurround\z@ \frozen@everymath{\mathgroup\macc@group\relax}%
    \macc@set@skewchar\relax
    \let\mathaccentV\macc@nested@a
    \if#31
    \macc@nested@a\relax111{#1}%
    \else
    \def\gobble@till@marker##1\endmarker{}%
    \futurelet\first@char\gobble@till@marker#1\endmarker
    \ifcat\noexpand\first@char A\else
    \def\first@char{}%
    \fi
    \macc@nested@a\relax111{\first@char}%
    \fi
    \endgroup
}
\makeatother

\title{A Wall Crossing Formula for %Quadratic Field Extensions}
Motivic Enumerative Invariants}
\author{Andrés Jaramillo Puentes}

\begin{document}

\begin{abstract}
    We prove an analog of the wall crossing formula for Wel\-schin\-ger invariants relating the difference of signed curve counting of real curves passing through configurations that differ by a pair of complex conjugated points, and a correspondence Welschinger invariant of the blow up.

    We prove this analogue for the motivic count of rational curves of fixed degree passing through a generic configuration of points, counted with a motivic multiplicity in the Grothendieck-Witt ring of a base field, extending the notions in the correspondence theorem between motivic invariants for $k$-rational point conditions and tropical curves, presented in~\cite{JPP23}.

    We use this formula to compute the degree 4 motivic enumerative invariants of the projective plane counting curves passing through configurations of points defined over quadratic extensions of a base field.
\end{abstract}
\maketitle
{\hypersetup{linkcolor=black}\tableofcontents}
\section{Introduction}

Let $k$ be a perfect field of characteristic different than $2$ and $3$. 
The Grothendieck-Witt ring $\GW(k)$  of the field $k$ is the group completion of the semi-ring of isometry classes of non-degenerate symmetric bilinear forms over $k$ under the direct sum, endowed with the tensor product. It is generated by the rank one forms
\[
\qinv{a}\colon k\times k\longrightarrow k \colon (x,y)\longmapsto axy, a\in k^{*}.
\]
%The Grothendieck-Witt ring is where our motivic enumerative invariants take place. For a natural number $n$, we denote by $n\cdot \qinv{a}$ the form with Gram matrix $a\mathrm{Id}_n$.

Given a natural number $d$ and a generic configuration $\mathcal{P}$ of $3d-1$ points of the projective plane over $k$, we define the motivic count
\[
N_{d,k}^{\Au}(\mathcal{P})\coloneqq\sum_{C} \mult(C)\in\GW(k),
\]
where the sum runs over rational degree $d$ curves $C$ in $\Pk$ passing through the configuration $\mathcal{P}$, counted with the multiplicity
\begin{equation}\label{eq:mult}
    \mult(C)\coloneqq\Tra[k(C)/k]{\prod_{z \text{ node of } C}\qinv{N_{k(z)/k(C) }D(z)}}\in\GW(k),
\end{equation}
where $k(C)$ is the field of definition of the curve $C$, $k(z)$ is the field of definition of the node and $D(z)\in k(z)^*/\left(k(z)^*\right)^2$ is the element such that the two tangent directions at the node $z$ define a quadratic field extension $k(z)\left[\sqrt{D(z)}\right]$ of $k(z)$. 
Here $\operatorname{Tr}$ denotes the map induced by the trace of the field extension on the Grothendieck-Witt rings and ${N}$ denotes the norm of the field extension. For details about this definition please confer to ~\cite[Subsection 1.2.4.]{KLSWCount}.

This definition can be naturally extended to other surfaces, generalizing the notion of degree by considering curves realizing a fixed cohomology class. 
In \cite{KLSWCount} and \cite{KLSWOrientation} is proved that the counts $N_{d,k}^{\Au}$ are independent of the choice of configuration of point as long as the configuration is generic; albeit in a more general statement.

For an algebraically closed field, since all finite field extensions are trivial, the configuration $\mathcal{P}$ is comprised of $k$-rational points and all motivic multiplicities are isometric to the symmetric bilinear form $\qinv{1}$.  In particular, over the complex numbers we have that
\[N_{d,\C}^{\Au}=N_d\cdot\qinv{1}\in\GW(\C),\]
where $N_d$ is the classical degree $d$ Gromov-Witten invariant of the complex projective plane. Therefore, the rank of the motivic Gromov-Witten invariant recovers the classical complex one.

Over the real numbers, the naïve count of curves depends on the configuration of points $\mathcal{P}$ and hence it is not an invariant. In \cite{Wel03}, it is proved that the invariance is restored by considering the real rational curves~$C$, endowed with multiplicities given by $\Wel_{\R}(C)\coloneqq(-1)^{e(C)}$, where $e(C)$ is the number of isolated real nodes of the curve~$C$. The signed sum of degree $d$ rational real curves passing through $\mathcal{P}$ is known as the degree $d$ Welschinger invariant~$W_d$ of the real projective plane. 
The motivic count considers curves defined over the base field as well as curves defined over finite field extensions. The total number of curves is the Gromov-Witten invariant $N_d$. However, there are two cases to consider. On one hand, a real rational curve $C$ has motivic multiplicity 
\[\mult[\R](C)=\qinv{\Wel_{\R}(C)}.\] 
On the other hand, a complex curve $C$ that is not defined over the real numbers has motivic multiplicity 
\[\mult[\R](C)=\Tra[\C/\R]{\qinv{1}}=\qinv{1}+\qinv{-1}\in\GW(\R).\]
In the former case, the signature of the motivic multiplicity coincides with the sign $\Wel_{\R}(C)$. In the the latter case, the signature of the motivic multiplicity vanishes. Additionally, let us remark that in this case the rank of the motivic multiplicity is two, accounting for the curve itself and for its complex conjugate. They are the same curve from a schematic point of view. Since the Welschinger invariant counts only real curves, we have that
\[\sig{N_{d,\R}^{\Au}(\mathcal{P})}=W_d.\]
By the Sylvester theorem, a non-degenerated real symmetric bilinear form is determined by its rank $N_d$ and its signature $W_d$. Hence,
\[N_{d,\R}^{\Au}(\mathcal{P})=
\frac{N_d+W_d}{2}\qinv{1}+\frac{N_d-W_d}{2}\qinv{-1}\in\GW(\R).\]

%Kontsevich recursion formula expresses the complex invariant $N_d$ in terms of the invariants of lower degree. 

Tropical geometry provides combinatorial tools for the computation of the complex and real invariants. We aim to harness them in order to determine relations among motivic invariants.
Mikhalkin’s celebrated correspondence theorem in \cite{Mi03} establishes a
correspondence between complex algebraic curves and tropical curves counted with 
multiplicities.
In this correspondence, we study the number of complex curves that tropicalize to a single tropical curve. This number is the complex multiplicity of the tropical curve. Subsequently, the complex invariant $N_d$ equals the sum over all tropical curves satisfying the required constrains (degree, genus and point conditions) counted with this complex multiplicity.
Additionally, this correspondence extends to the real invariants. One can study the local contribution to the real signed count, defining a real multiplicity in such a way that the real invariant $W_d$ equals the sum over the same set of tropical curves counted with the real multiplicity.
In \cite{JPP23}, we proved Theorem~\ref{thm:motcor}, extending the correspondence between algebraic and tropical curves to arbitrary perfect fields of large characteristic. In order to state this theorem, let us first introduce  the notion of motivic multiplicity for tropical curves.

Assume $\Gamma$ is a simple degree $\Delta$ embedded tropical curve, i.e., its dual subdivision  $S$ is comprised of triangles and parallelograms.
Additionally, assume that all unbound edges of $\Gamma$ have weight one.

\begin{df}
    Let $v$ be a trivalent vertex of $\Gamma$ dual to the triangle $\Delta_v$. Let us denote by $\sigma_v$, $\sigma_v'$ and $\sigma_v''$ the three edges of the triangle $\Delta_v$ and by~$\operatorname{Int}(\Delta_v)$ its number of interior lattice points. % of $\Delta_v$.
    The \emph{motivic multiplicity} of the vertex~$v$ is the class
\[m_v^{\A^1}\coloneqq\begin{cases}
    \qinv{ (-1)^{\operatorname{Int}(\Delta_v)}\vert\sigma_v\vert\vert\sigma'_v\vert\vert\sigma''_v\vert}+\frac{\vert\Delta_v\vert-1}{2}h&\text{if $\Delta_v$ has only odd edges,}\\
    \frac{\vert\Delta_v\vert}{2}h&\text{otherwise.}
\end{cases}\]
Here,  the norm $\vert\cdot\vert$ denotes the lattice length for a segment, and twice the Euclidean area for the triangle $\Delta_v$, and $h$ denotes the hyperbolic form $\qinv{1}+\qinv{-1}\in\GW(k)$.
The \emph{motivic multiplicity} of the tropical curve $\Gamma$ is the product of its vertex multiplicities over all its trivalent vertices
    \[\mult(\Gamma)\coloneqq\prod_{v\in V(\Gamma), \operatorname{val}(v)=3}\mkern-15mu m_v^{\A^1}.\]
\end{df}

Let us remark that the rank and signature (for real fields) of the motivic  multiplicity is the complex and real multiplicity, respectively, both for vertices and, subsequently, for simple nodal tropical curves.

Now, in order to study the relation between algebraic and tropical curves, we consider our algebraic curves as curves in families. We see the coefficients of a defining equation of the curve as elements of the Puiseux series with coefficients over $k$. The inclusion of fields $k\subset k\Puiseux$ induces an isomorphism of rings (cf. \cite[Theorem 4.7]{MPS22})
\[
\begin{array}{ccc}
\GW(k)&\longrightarrow&\GW(k\Puiseux)\\
\qinv{a}&\longmapsto&\qinv{a}
\end{array}
\]
which allow us to compute the motivic invariants we are studying in the Puiseux series and mapping back to the base field through this canonical isomorphism.
Let us consider
 a generic configuration ${\mathcal{P}}$ of $k$-rational points in $\mathbb{P}^2_k$,
 a generic configuration of $k\Puiseux$-rational points $\widetilde{\mathcal{P}}$ in $\mathbb{P}^2_{k\Puiseux}$ such that the set of initial terms %$\operatorname{it}\left(\widetilde{\mathcal{P}}\right)\coloneqq 
$\{a_{i_0}\mid \sum_{i=i_0}^\infty a_it^{i/N}\in \mathcal{P} \}$ of $\widetilde{\mathcal{P}}$ coincides with $\mathcal{P}$ and such that the set of point-wise valuations 
$\overline{\mathcal{P}}\coloneqq\val\left(\widetilde{\mathcal{P}}\right)$ defined as the set 
\[\displaystyle\left\{(-i_0/N,-j_0/M) \middle| \left(\sum_{i=i_0}^\infty a_it^{i/N},\sum_{j=j_0}^\infty b_jt^{j/M}\right) \in \widetilde{\mathcal{P}},a_{i_0},b_{j_0}\neq0 \right\}
\]
is tropically generic. This context allows us to state the following correspondence theorem.
\begin{thm}[\cite{JPP23}]\label{thm:motcor}
    Let $k$ be a field of characteristic $0$ or characteristic greater than $d$.
If $\Gamma\subset\R^2$ is a rational degree $d$ tropical curve passing through $\overline{\mathcal{P}}$, then under the canonical isomorphism~$\GW(k\Puiseux)\cong \GW(k)$ the quadratic multiplicity $\mult(\Gamma)$ of $\Gamma$ is mapped to
    \[\sum_C \mult(C),\]
    where the sum runs over all rational curves $C$ in $\mathbb{P}^2_{k\Puiseux}$ passing through $\widetilde{\mathcal{P}}$ and tropicalizing to $\Gamma$.
\end{thm}

A consequence of Theorem~\ref{thm:motcor} together with the invariance of the motivic counts with respect to the configuration of points, proved in~\cite{KLSWCount}, yields the following corollary.

\begin{coro}[\cite{JPP23}] \label{coro:cthm}
    If $k$ is a perfect field of characteristic $0$ or characteristic greater than $d$, then 
\[N_{d,k}^{\Au}(\mathcal{P})=
\frac{N_d+W_d}{2}\qinv{1}+\frac{N_d-W_d}{2}\qinv{-1}\in\GW(k),\]
for any generic configuration $\mathcal{P}$ of $3d-1$ $k$-rational points in $\mathbb{P}^2_k$.
\end{coro}
In motivic enumerative geometry, 
whenever the enumerative problem can be stated over $\Z$ in a smooth and proper manner the invariants lie in the ideal generated by the forms $\qinv{1}$ and $\qinv{-1}$, since they generate the Grothendieck-Witt ring of the integers. 
However, this is not always the case. For example, the count of lines on the Fermat cubic surface yields $\qinv{-3} + 4h + \qinv{1}+ \operatorname{Tr}_{k(\sigma)/k}{\qinv{1}}$, where $\operatorname{Tr}_{k(\sigma)/k}$ is the trace form from the étale algebra given by the fields of definition of the point conditions $\sigma$ to the base field. 
Another example is the count of rational plane cubics
\[N_{3,k}^{\Au}(\mathcal{P})=
2h+\operatorname{Tr}_{k(\sigma)/k}{\qinv{1}},\]
for any generic configuration $\mathcal{P}$ (cf. \cite{KLSWCount} for details about the two examples).
If the configuration of points $\mathcal{P}$ is entirely $k$-rational, then the field of definition of each point condition is $k$ and the étale algebra $k(\sigma)=k^8$; and hence $\operatorname{Tr}_{k(\sigma)/k}{\qinv{1}}=8\qinv{1}$, which coincides with the statement in Corollary~\ref{coro:cthm}.

The goal of this article is to extend our computation of $N_{d,k}^{\Au}$ to quadratic extensions using wall crossing formulas. 
I.e., we aim to compute $N_{d,k}^{\Au}$ for generic configurations of points $\mathcal{P}_{\Bar{c}}$, where $\Bar{c}=(c_i)_{i\in I}$, comprised of a mix of~$k$-rational points~$p_i$ and points~$q_j$ defined over quadratic extensions $k(\sqrt{c_j})$ of~$k$. The point~$q_j$ is a scheme-theoretical point, it represents two points over the algebraic closure of $k$.

All examples presented in \cite{KLSWCount} are either constant or linear in the traces over the field of definition of the point conditions. Theorem~\ref{thm:WCF} implies that~$N_{d,k}^{\Au}$ is a monic polynomial of degree $(d-1)(d-2)/2$, as exemplified in Section~\ref{sec:apps}. Or, in a more general setting, the count $N_{X,D,\sigma}$ of the degree~$D$ rational curves through points defined over quadratic field extensions $\sigma$ in a toric del Pezzo surface $X$ is a monic polynomial in the traces over the field of definition of the point conditions, of degree $(D^2-c_1(X)D+2)/2$.

Over the real numbers, one can partition the configuration of points. Pick natural numbers $s$ and $t$ such that~$t+2s=3d-1$. Then, for a generic configuration $\mathcal{P}_s$ comprised of $t$ real points in $\mathbb{P}^2_{\R}$ and $s$ pairs of complex conjugated points in $\mathbb{P}^2_{\C}\setminus \mathbb{P}^2_{\R}$, we can define $W_d(s)$ as the count of real curves passing through $\mathcal{P}_s$ counted with the sign $\Wel_{\R}$.
In~\cite{Wel03}, Welschinger proved that the difference of real invariants differing by a pair of complex conjugated points in their point conditions, equals  twice the real invariant associated to the blow up, i.e., they satisfy
\begin{equation}\label{eq:WCr}
  W_d(s)-W_d(s+1)=2\cdot\widetilde{W}_{d-2E}(s),  
\end{equation}
where $\widetilde{W}_{d-2E}(s)$ is the signed count of real curves $C$ lying on the blow up~$\operatorname{Bl}_p(\mathbb{P}^2_{\R})$ of~$\mathbb{P}^2_{\R}$ at a real point $p\in \mathcal{P}_s$ passing through a generic configuration of~$t-2$ real points and $s$ pairs of complex conjugated points, such that the cohomology class of $C$ equals $dL-2E$, where $L$ is the strict transform of a generic line in $\mathbb{P}^2_{\R}$ and $E$ is the exceptional divisor of the blow up.

Over arbitrary fields, we study the difference of motivic invariants associated to configurations that differ by one pair of $k$-rational points on one hand, and a pair of conjugated points in a quadratic extension on the other. Over the algebraic closure, there is no distinction between the two cases and so, the rank of such difference should vanish. Over the real number such difference must specialise to the Equation~\eqref{eq:WCr}.

In order to state our results, let us extend our definitions. 
Let $X$ be a toric del Pezzo surface such that the blow up of $X$ at a $k$-point~$p\in X$ is also a toric del Pezzo surface, together with an appropriated notion of degree. Hence, assume that $X$ is the projective plane or its blow up at one or two $k$-points. Let $D=dL-\sum d_{E_i}E_i$ be a divisor in $X$ such that~$n\coloneqq\deg(-D\cdot K_{X})-1\geq0$. Here $d,d_{E_i}\in\Z$, $L$ is the strict transform of a generic line in $\mathbb{P}^2_k$, and $E_i$ is the exceptional divisor of the blow up, if any.
Given a generic configuration~$\mathcal{P}$ of $n$ points of $X$, we define the motivic count 
\[
N_{X,D,k}^{\Au}(\mathcal{P})\coloneqq\sum_{C} \mult(C)\in\GW(k),
\]
where the sum runs over all rational curves $C$ in $X$, in the same cohomology class as $D$, passing through $\mathcal{P}$, counted with the multiplicity given by Equation~\eqref{eq:mult}. By \cite{KLSWCount,KLSWOrientation} this is invariant with respect to the generic configuration as long as the points are defined over isomorphic field extensions.
If two configurations differ by a pair of points, either $k$-points or a pair of conjugated points, we have the following relation.

\begin{thm}\label{thm:WCF}
    If  $\mathcal{P}$ and $\mathcal{P'}$ are two generic configurations of points of~$X$ such that $\mathcal{P}\setminus\mathcal{P'}$ consist of two $k$-rational points $p, p'$, and $\mathcal{P}'\setminus\mathcal{P}$ consist of a pair of conjugated points $q$ defined over a quadratic extension $k(\sqrt{c})$, then
\[N^{\mathbb{A}^1}_{X,D,k}(\mathcal{P})-N^{\mathbb{A}^1}_{X,D,k}(\mathcal{P}')=
\delta \cdot N^{\mathbb{A}^1}_{\operatorname{Bl}_p(X),D-2E,k}(\mathcal{P}\cap \mathcal{P}'),\]
where 
%\[\beta:=\Tra[k(\sqrt{c}))/k]{\qinv{1}},\]
$\delta=%2\qinv{1}-\beta=
\Tra[k\times k/k]{\qinv{1}}-\Tra[k(\sqrt{c})/k]{\qinv{1}}$, and $E$ is the exceptional divisor of the blow up.
\end{thm}

The rank of the form $\delta$ is zero and its signature, for real fields and negative~$c$, equals $2$.
Let us remark that the rank of both sides of the formula vanishes, and taking the signature (for real fields and negative $c$) specialise our statement to Welschinger's theorem in Equation~\eqref{eq:WCr}.

Let us comment that although Theorem~\ref{thm:WCF} is stated for algebraic invariants, we can use the same arguments for the tropical invariants $N^{\mathbb{A}^1,\text{trop}}$ studied in~\cite{JPMPR}. 
In a joint project with Markwig, Pauli and  Röhrle, we are working towards a correspondence theorem for quadratic field extensions, from which this formula can be deduced.
Brugallé and Wickelgren are working towards a quadratically enriched Abramovich-Bertram formula that would generalise this relation.

In Section~\ref{sec:WCF} we provide some key facts of tropical, toric and $\mathbb{A}^1$-enu\-mer\-ative geometry. Followed by a proof of Theorem~\ref{thm:WCF}.
In Section~\ref{sec:apps} we  compute the motivic invariants of degree 4 of $\mathbb{P}^2_k$ for configurations of points having points defined over quadratic field extensions.

\subsection{Acknowledgements}
The author was supported by the ERC programme QUADAG.  This paper is part of a project that has received funding from the European Research Council (ERC) under the European Union's Horizon 2020 research and innovation programme (grant agreement No. 832833). The author would like to thank the Universität Duisburg-Essen and the Università degli Studi di Napoli Federico II for support.\\ 
\includegraphics[width=2.3cm]{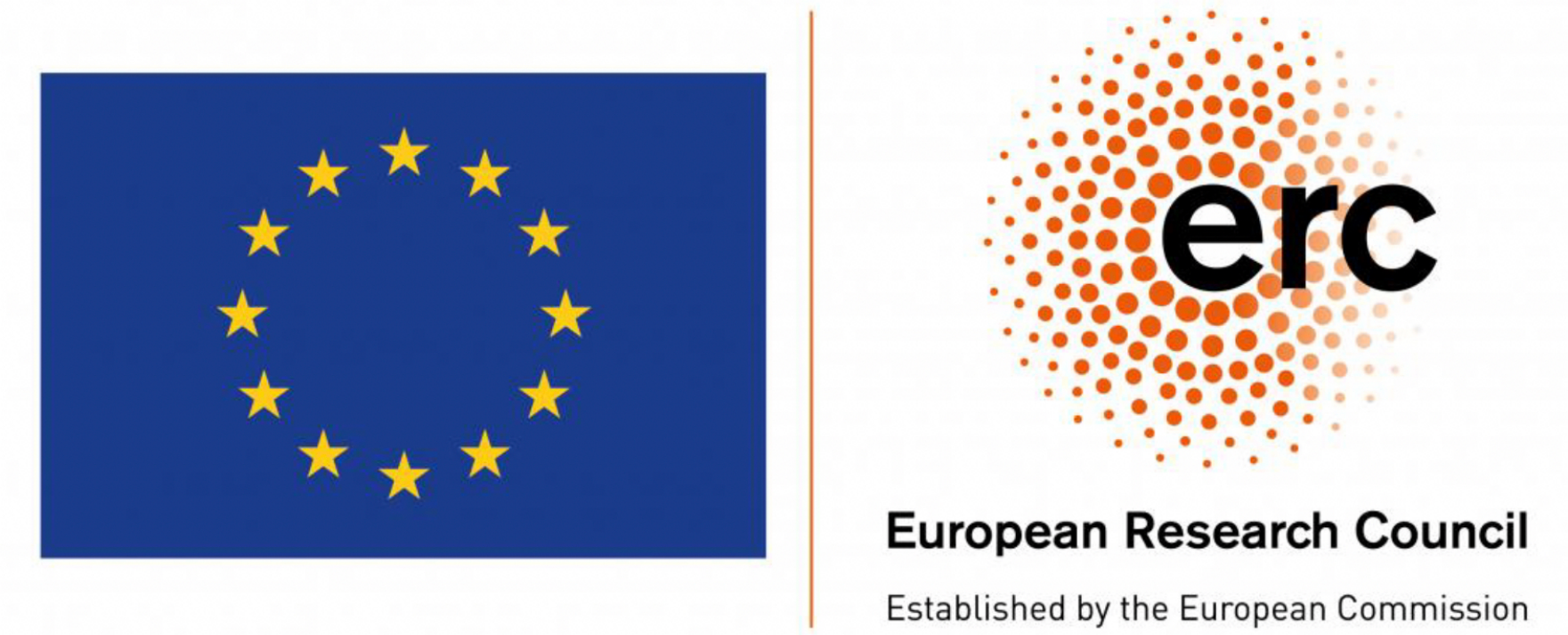}

The author would like to thank Erwan Brugallé, Hannah Markwig, Marc Levine, Sabrina Pauli and Kirsten Wickelgren, for fruitful discussions about the different techniques, motivic and tropical, that were used in this paper.

\section{Wall crossing formula for quadratic field extensions}
\label{sec:WCF}

In the first part of this section we provide two notions of tropical geometry and some relations between elements of the Grothendieck-Witt ring that are key to the proof of Theorem~\ref{thm:WCF} presented at the end of this section.

%\section{Technical lemmas}
\begin{figure}
\centering
\definecolor{xdxdff}{rgb}{0.6,0.6,1.0}
\definecolor{ffffqq}{rgb}{1.0,1.0,0.0}
\definecolor{qqqqff}{rgb}{0.0,0.0,1.0}
\definecolor{zzttqq}{rgb}{0.6,0.2,0.0}

\begin{tikzpicture}[line cap=round,line join=round,>=triangle 45,x=1.0cm,y=1.0cm,scale=0.8]

\clip(-.5,1.3) rectangle (11,5);
\fill[line width=2.4pt,color=zzttqq,fill=zzttqq,fill opacity=0.10000000149011612] (0.,0.) -- (4.,0.) -- (0.,4.) -- cycle;
\fill[line width=2.4pt,color=zzttqq,fill=zzttqq,fill opacity=0.10000000149011612] (6.,0.) -- (6.,2.) -- (8.,2.) -- (10.,0.) -- cycle;
\fill[line width=2.4pt,color=zzttqq,fill=zzttqq,fill opacity=0.10000000149011612] (12.,2.) -- (14.,2.) -- (14.,0.) -- (12.,0.) -- cycle;
\draw [line width=2.4pt,color=zzttqq] (0.,0.)-- (4.,0.);
\draw [line width=2.4pt,color=zzttqq] (4.,0.)-- (0.,4.);
\draw [line width=2.4pt,color=zzttqq] (0.,4.)-- (0.,0.);
\draw [line width=2.4pt,color=zzttqq] (6.,0.)-- (6.,2.);
\draw [line width=2.4pt,color=zzttqq] (6.,2.)-- (8.,2.);
\draw [line width=2.4pt,color=zzttqq] (8.,2.)-- (10.,0.);
\draw [line width=2.4pt,color=zzttqq] (10.,0.)-- (6.,0.);
\draw [line width=2.4pt,color=zzttqq] (12.,2.)-- (14.,2.);
\draw [line width=2.4pt,color=zzttqq] (14.,2.)-- (14.,0.);
\draw [line width=2.4pt,color=zzttqq] (14.,0.)-- (12.,0.);
\draw [line width=2.4pt,color=zzttqq] (12.,0.)-- (12.,2.);

\draw [fill=xdxdff!50] (6.,3.) circle (2.5pt);
\draw [fill=xdxdff!50] (6.,4.) circle (2.5pt);
\draw [fill=xdxdff!50] (7.,3.) circle (2.5pt);

\draw [fill=xdxdff] (0.,0.) circle (2.5pt);
\draw [fill=xdxdff] (4.,0.) circle (2.5pt);
\draw [fill=xdxdff] (0.,4.) circle (2.5pt);
\draw [fill=xdxdff] (1.,0.) circle (2.5pt);
\draw [fill=xdxdff] (2.,0.) circle (2.5pt);
\draw [fill=xdxdff] (3.,0.) circle (2.5pt);
\draw [fill=xdxdff] (3.,1.) circle (2.5pt);
\draw [fill=xdxdff] (2.,1.) circle (2.5pt);
\draw [fill=xdxdff] (1.,1.) circle (2.5pt);
\draw [fill=xdxdff] (0.,1.) circle (2.5pt);
\draw [fill=xdxdff] (0.,2.) circle (2.5pt);
\draw [fill=xdxdff] (1.,2.) circle (2.5pt);
\draw [fill=xdxdff] (2.,2.) circle (2.5pt);
\draw [fill=xdxdff] (1.,3.) circle (2.5pt);
\draw [fill=xdxdff] (0.,3.) circle (2.5pt);
\draw [fill=xdxdff] (6.,0.) circle (2.5pt);
\draw [fill=xdxdff] (6.,2.) circle (2.5pt);
\draw [fill=xdxdff] (8.,2.) circle (2.5pt);
\draw [fill=xdxdff] (10.,0.) circle (2.5pt);
\draw [fill=xdxdff] (12.,2.) circle (2.5pt);
\draw [fill=xdxdff] (14.,2.) circle (2.5pt);
\draw [fill=xdxdff] (14.,0.) circle (2.5pt);
\draw [fill=xdxdff] (12.,0.) circle (2.5pt);
\draw [fill=xdxdff] (7.,2.) circle (2.5pt);
\draw [fill=xdxdff] (6.,1.) circle (2.5pt);
\draw [fill=xdxdff] (7.,1.) circle (2.5pt);
\draw [fill=xdxdff] (8.,1.) circle (2.5pt);
\draw [fill=xdxdff] (7.,0.) circle (2.5pt);
\draw [fill=xdxdff] (8.,0.) circle (2.5pt);
\draw [fill=xdxdff] (9.,0.) circle (2.5pt);
\draw [fill=xdxdff] (9.,1.) circle (2.5pt);
\draw [fill=xdxdff] (12.,1.) circle (2.5pt);
\draw [fill=xdxdff] (13.,1.) circle (2.5pt);
\draw [fill=xdxdff] (14.,1.) circle (2.5pt);
\draw [fill=xdxdff] (13.,2.) circle (2.5pt);
\draw [fill=xdxdff] (13.,0.) circle (2.5pt);

\draw (8.,-1.) node{$\Sigma_1$};
\draw (2.,-1.) node{$\mathbb{P}^2$};
\draw (13.,-1.) node{$\operatorname{Bl}_p\left(\Sigma_1\right)$};
\end{tikzpicture}
    \caption{Lattice polygons corresponding to the blow up at a $k$-point fixed by the torus action.}
    \label{fig:polysbu}
\end{figure}
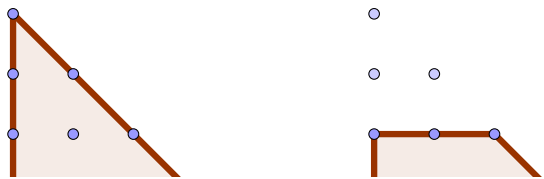

\subsection{Tropical Geometry}
\label{sec:trop}
This paper assumes a basic knowledge of tropical geometry, like tropical curves, their degree, genus and the dual subdivision of its degree. We refer the reader to~\cite[Section 3]{JPP22}, where we give an introduction to the basic notions of tropical geometry that we use in this paper, or to more classical references like~\cite{Mi06,IMS09} for further details.
In this subsection we introduce two notions that we use in our arguments that may be unfamiliar to the reader.

The first one is the notion of \emph{vertically stretched} configurations of points in $\R^2$. These are configurations whose points are ordered, with monotone entries in both its abscissae and ordinates, and such that the difference of consecutive ordinates is of larger magnitude than the difference of the corresponding abscissae. Tropical curves that pass through vertically stretched configurations of points in $\R^2$ were studied in~\cite{FM09}. 
Their combinatorics can be studied systematically. Vertically stretched configurations form an open set in the space of configurations, which allows us to pick a configuration of this sort and still study generic conditions on tropical curves. 

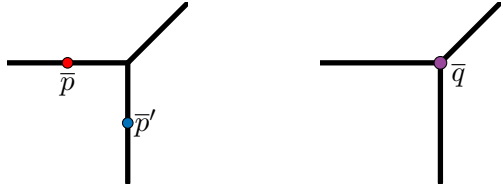
\begin{figure}
\centering
\begin{tikzpicture}[line cap=round,line join=round,>=triangle 45,x=1.0cm,y=1.0cm, scale=0.8]
\clip(0,0) rectangle (3,3);
\draw [line width=2.pt] (2.,0.)-- (2.,2.);
\draw [line width=2.pt] (2.,2.)-- (0.,2.);
\draw [line width=2.pt] (2.,2.)-- (3.,3.);
\draw [fill=red] (1.,2.) circle (2.5pt);
\draw (1,1.65) node {$\overline{p}$};
\draw [fill=RoyalBlue] (2.,1.) circle (2.5pt);
\draw (2.3,1) node {$\overline{p}'$};
\end{tikzpicture}
\hspace{1.5cm}
\begin{tikzpicture}[line cap=round,line join=round,>=triangle 45,x=1.0cm,y=1.0cm,scale=0.8]
\clip(0,0) rectangle (3,3);
\draw [line width=2.pt] (2.,0.)-- (2.,2.);
\draw [line width=2.pt] (2.,2.)-- (0.,2.);
\draw [line width=2.pt] (2.,2.)-- (3.,3.);
\draw [fill=Purple] (2.,2.) circle (3pt);
\draw (2.3,1.8) node {$\overline{q}$};
\end{tikzpicture}
    \caption{Tropical curves of degree $\Delta_1$ passing through two $k$-points or a double point defined over a quadratic extension $k(\sqrt{c})$.}
    \label{fig:alpha}
\end{figure}

A second fact that we use is that if a tropical curve having degree~$\Delta$ as in the left polygon in Figure~\ref{fig:polysbu} passes through a vertically stretched configuration of points such that two simple points $\overline{p},\overline{p}'$ or a \emph{double} point (a point arising from the degeneration of two simple points) $\overline{q}'$ that are at the topmost part of the configuration, then the tropical curve around these points looks like the ones in Figures~\ref{fig:alpha}, \ref{fig:beta} and \ref{fig:beta2}.

\begin{figure}
\centering
\definecolor{xdxdff}{rgb}{0.6,0.6,1.0}
\definecolor{ffffqq}{rgb}{1.0,1.0,0.0}
\definecolor{qqqqff}{rgb}{0.0,0.0,1.0}
\definecolor{zzttqq}{rgb}{0.6,0.2,0.0}

\begin{tikzpicture}[line cap=round,line join=round,>=triangle 45,x=1.0cm,y=1.0cm,scale=0.8]

\clip(-.5,-2) rectangle (11,5);
\fill[line width=2.4pt,color=zzttqq,fill=zzttqq,fill opacity=0.10000000149011612] (0.,0.) -- (4.,0.) -- (0.,4.) -- cycle;
\fill[line width=2.4pt,color=zzttqq,fill=zzttqq,fill opacity=0.10000000149011612] (6.,0.) -- (6.,2.) -- (8.,2.) -- (10.,0.) -- cycle;
\fill[line width=2.4pt,color=zzttqq,fill=zzttqq,fill opacity=0.10000000149011612] (12.,2.) -- (14.,2.) -- (14.,0.) -- (12.,0.) -- cycle;
\draw [line width=2.4pt,color=zzttqq] (0.,0.)-- (4.,0.);
\draw [line width=2.4pt,color=zzttqq] (4.,0.)-- (0.,4.);
\draw [line width=2.4pt,color=zzttqq] (0.,4.)-- (0.,0.);
\draw [line width=2.4pt,color=zzttqq] (6.,0.)-- (6.,2.);
\draw [line width=2.4pt,color=zzttqq] (6.,2.)-- (8.,2.);
\draw [line width=2.4pt,color=zzttqq] (8.,2.)-- (10.,0.);
\draw [line width=2.4pt,color=zzttqq] (10.,0.)-- (6.,0.);
\draw [line width=2.4pt,color=zzttqq] (12.,2.)-- (14.,2.);
\draw [line width=2.4pt,color=zzttqq] (14.,2.)-- (14.,0.);
\draw [line width=2.4pt,color=zzttqq] (14.,0.)-- (12.,0.);
\draw [line width=2.4pt,color=zzttqq] (12.,0.)-- (12.,2.);

\draw [fill=xdxdff] (0.,0.) circle (2.5pt);
\draw [fill=xdxdff] (4.,0.) circle (2.5pt);
\draw [fill=xdxdff] (0.,4.) circle (2.5pt);
\draw [fill=xdxdff] (1.,0.) circle (2.5pt);
\draw [fill=xdxdff] (2.,0.) circle (2.5pt);
\draw [fill=xdxdff] (3.,0.) circle (2.5pt);
\draw [fill=xdxdff] (3.,1.) circle (2.5pt);
\draw [fill=xdxdff] (2.,1.) circle (2.5pt);
\draw [fill=xdxdff] (1.,1.) circle (2.5pt);
\draw [fill=xdxdff] (0.,1.) circle (2.5pt);
\draw [fill=xdxdff] (0.,2.) circle (2.5pt);
\draw [fill=xdxdff] (1.,2.) circle (2.5pt);
\draw [fill=xdxdff] (2.,2.) circle (2.5pt);
\draw [fill=xdxdff] (1.,3.) circle (2.5pt);
\draw [fill=xdxdff] (0.,3.) circle (2.5pt);
\draw [fill=xdxdff] (6.,0.) circle (2.5pt);
\draw [fill=xdxdff] (6.,2.) circle (2.5pt);
\draw [fill=xdxdff] (8.,2.) circle (2.5pt);
\draw [fill=xdxdff] (10.,0.) circle (2.5pt);
\draw [fill=xdxdff] (12.,2.) circle (2.5pt);
\draw [fill=xdxdff] (14.,2.) circle (2.5pt);
\draw [fill=xdxdff] (14.,0.) circle (2.5pt);
\draw [fill=xdxdff] (12.,0.) circle (2.5pt);
\draw [fill=xdxdff] (7.,2.) circle (2.5pt);
\draw [fill=xdxdff] (6.,1.) circle (2.5pt);
\draw [fill=xdxdff] (7.,1.) circle (2.5pt);
\draw [fill=xdxdff] (8.,1.) circle (2.5pt);
\draw [fill=xdxdff] (7.,0.) circle (2.5pt);
\draw [fill=xdxdff] (8.,0.) circle (2.5pt);
\draw [fill=xdxdff] (9.,0.) circle (2.5pt);
\draw [fill=xdxdff] (9.,1.) circle (2.5pt);
\draw [fill=xdxdff] (12.,1.) circle (2.5pt);
\draw [fill=xdxdff] (13.,1.) circle (2.5pt);
\draw [fill=xdxdff] (14.,1.) circle (2.5pt);
\draw [fill=xdxdff] (13.,2.) circle (2.5pt);
\draw [fill=xdxdff] (13.,0.) circle (2.5pt);

\draw (8.,-1.) node{$\Sigma_1$};
\draw (2.,-1.) node{$\mathbb{P}^2$};
\draw (13.,-1.) node{$\operatorname{Bl}_p\left(\Sigma_1\right)$};
\end{tikzpicture}
    \caption{Lattice polygons corresponding to quartic curves in the projective plane $\mathbb{P}^2$, and to curves in the first Hirzebruch surface $\Sigma_1$ of de\-gree~$4-2E$,  respectively.}
    \label{fig:polys}
\end{figure}
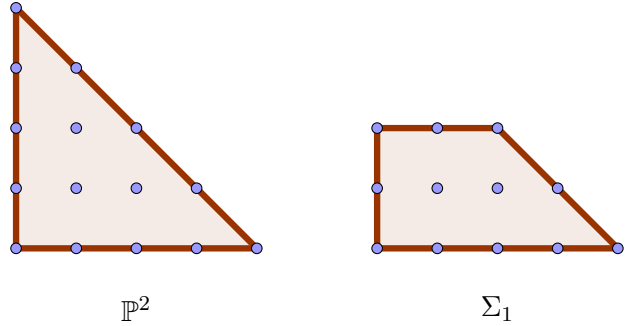

The second notion that we need is not proper to tropical geometry, but to toric geometry. 
If we assume that a toric surface $X$ has a zero-dimensional toric orbit at a point $p$, then the blow up $\operatorname{Bl}_p(X)$ is toric and their associated polygons are like the ones in Figure~\ref{fig:polysbu}. We refer the reader to~\cite{Fulton} for a detailed construction. 
Since we only consider the projective plane $\mathbb{P}^2$ and its blow up in up to three $k$-points, and we can always change coordinates to set any point as a point invariant by the toric action, the polygons of the blow ups look like the standard triangle $\Delta_m\coloneqq\operatorname{Conv}\{(0,0),(0,m),(m,0)\}$, to which we chop off up to three corners, each corner being isomorphic to $\Delta_2$ by a $SL(2,\Z)$ action.
Examples can be found in Figures~\ref{fig:polys} and~\ref{fig:polys2}.

\subsection{\texorpdfstring{$\mathbb{A}^1$-}{A1} Enumerative Geometry}

Motivic invariants are elements in the Gro\-then\-dieck-Witt ring $\GW(k)$ of the base field $k$, which we assumed perfect of characteristic not $2$ or $3$.
In this paper, we make use of the arithmetic of elements of $\GW(k)$ as well as the trace map from quadratic extensions to the base field \[\operatorname{Tr}_{k(\sqrt{c})/k}\colon\GW(k(\sqrt{c}))\longrightarrow\GW(k)\] induced by the algebraic trace $\operatorname{tr}_{k(\sqrt{c})/k}\colon k(\sqrt{c})\longrightarrow k$ that sends an element to the sum of its Galois conjugates. 
In this section we present the relations that we use and the aforementioned trace map. 
For further details about the notions of $\mathbb{A}^1$-enumerative geometry that we use in motivic tropical geometry, we refer the reader to~\cite[Section 2]{JPP22}.
As we stated in the introduction, the Grothendieck-Witt ring $\GW(k)$ is generated by the classes of bilinear forms 
\[\qinv{ a} \colon k\times k\rightarrow k\colon (x,y)\longmapsto axy,\]
for $a\in k^{*}$. These generators are subject to the relations
\begin{itemize}

\item $\qinv{ab^2}=\qinv{a}$,
    \item $\qinv{ a} \qinv{ b}=\qinv{ ab}$,
    \item $\qinv{ a}+\qinv{ b}=\qinv{ a+b}+\qinv{ ab(a+b)}$,
    \item $\qinv{a}+\qinv{-a}=\qinv{1}+\qinv{-1}$,
\end{itemize}
given that $a,b,a+b\in k^{*}$. Recall that we define they hyperbolic form $h$ as the form~$ \qinv{1}+\qinv{-1}$.
The following proposition computes the trace form of elements in $\GW(k(\sqrt{c}))$ that are multiples of one element of the base~$\{1,\sqrt{c}\}$ of the $k$-vector space~$k(\sqrt{c})$.

\begin{prop}[\cite{JPP22}]
If $a\in k$, then
\begin{enumerate}
    \item $\operatorname{Tr}_{k(\sqrt{c})/k}(\qinv{a})=\qinv{2a}+\qinv{2ac}$,
    \item $\operatorname{Tr}_{k(\sqrt{c})/k}(\langle{a\sqrt{c}}\rangle)=h$.
\end{enumerate}
\end{prop}

\begin{lm}\label{lm:betamult}
If $\alpha\in\GW(k(\sqrt{c}))$ is an integer sum of elements $\qinv{a}$, $a\in k^*$, then
\[\operatorname{Tr}_{k(\sqrt{c})/k}(\alpha)=\alpha\cdot\operatorname{Tr}_{k(\sqrt{c})/k}(\qinv{1})\in\GW(k).\]
\end{lm}

\begin{proof}
    The trace form map $\operatorname{Tr}_{k(\sqrt{c})/k}$ is a homomorphism of groups. Hence, the lemma follows from the equality
    \[
    \operatorname{Tr}_{k(\sqrt{c})/k}(\qinv{a})=\qinv{2a}+\qinv{2ac}=\qinv{a}\cdot\left(\qinv{2}+\qinv{2c}\right)
    =\qinv{a}\cdot\operatorname{Tr}_{k(\sqrt{c})/k}(\qinv{1}).
    \]
\end{proof}

In general, motivic multiplicities can be arbitrary. 
However, the motivic multiplicities of curves passing through configurations of points that are~$k$-rational and the curves under consideration in this paper are generated by elements $\qinv{a}, a\in k^*$. Thus, we have that in this setting the trace forms of quadratic extensions are multiplicative.

% \begin{thm}
%     If  $\mathcal{P}$ and $\mathcal{P'}$ are two generic configurations of points of~$X$ such that $\mathcal{P}\setminus\mathcal{P'}$ consist of two $k$-rational points $p, p'$, and $\mathcal{P}'\setminus\mathcal{P}$ consist of a pair of conjugated points $q$ defined over a quadratic extension $k(\sqrt{c})$, then
% \[N^{\mathbb{A}^1}_{X,D,k}(\mathcal{P})-N^{\mathbb{A}^1}_{X,D,k}(\mathcal{P}')=
% \delta \cdot N^{\mathbb{A}^1}_{\operatorname{Bl}_p(X),D-2E,k}(\mathcal{P}\cap \mathcal{P}'),\]
% where 
% %\[\beta:=\Tra[k(\sqrt{c}))/k]{\qinv{1}},\]
% \[\delta=%2\qinv{1}-\beta=
% \Tra[k\times k/k]{\qinv{1}}-\Tra[k(\sqrt{c}))/k]{\qinv{1}}.\]
% \end{thm}

\subsection{Proof of Theorem~\ref{thm:WCF}}
In \cite{KLSWCount}, it is proved that the motivic count $N^{\mathbb{A}^1}_{X,D,k}$ does not depend on the generic configuration, but only on the isomorphism class of the étale algebra of fields of definition of the point conditions. Hence, without lost of generality, we can pick our configurations~$\mathcal{P}$ and~$\mathcal{P}'$ in such a way that the respective configurations over the field of Puiseux series $\widetilde{P}$ and $\widetilde{P}'$ are mapped by the valuation to a configuration of points that is vertically stretched (see preamble of Theorem~\ref{thm:motcor}) and that send, under the valuation map, the symmetric difference of the two configurations to the topmost part. Explicitly, we want to pick configurations $\widetilde{P}$ and $\widetilde{P}'$ such that the points~$\overline{p}\coloneqq \operatorname{val}(\widetilde{p}),\overline{p}'\coloneqq \operatorname{val}(\widetilde{p}')$ and $\overline{q}\coloneqq\operatorname{val}(\widetilde{q})$ are the topmost part of a vertically stretched configuration, where $\widetilde{p},\widetilde{p}'$ and $\widetilde{q}$ are lifts of $p, p'\in\mathcal{P}\setminus\mathcal{P}'$ and~$q\in\mathcal{P}'\setminus\mathcal{P}$ to~$k\Puiseux$ and~$k(\sqrt{c})\Puiseux$, respectively.

Since we picked our configuration so that the points under study correspond to the corner of the polygon that is being blown up,
this choice of configuration constrains the tropical curves that pass through the configurations $\overline{\mathcal{P}}$ and $\overline{\mathcal{P}}'$ to have a dual polygon divided into two different combinatorial patterns: we can chop off a  triangle $\Delta_1$ or $\Delta_2$ in such a way that the curve corresponding to the triangle $\Delta_i$ contains $\overline{p},\overline{p}'$ or $\overline{q}$, and the remaining curve contains $\overline{\mathcal{P}}\cap\overline{\mathcal{P}}'$.

The first case consists of tropical curves that pass through the special points as in Figure~\ref{fig:alpha}.
Cutting the plane $\R^2$ at a horizontal line just below the special points will split each tropical curve $\Gamma$ into two pieces: the top part $\Gamma^{\text{t}}$ and the bottom part $\Gamma_{\text{b}}$. The curve $\Gamma^{\text{t}}$ has degree~$\Delta_1$, while the curve~$\Gamma_{\text{b}}$ has as degree the Newton polygon $\overline{\Delta}\coloneqq\Delta\setminus\Delta_1$. 
The algebraic curves that tropicalize to $\Gamma$ split into pieces lying on $\operatorname{Tor}(\Delta_v)$ for every vertex~$v$ of $\Gamma$. In $\operatorname{Tor}(\Delta_1)$, the corresponding piece is a line that tropicalizes to~$\Gamma^{\text{t}}$. This line has no nodes, it is defined over the base field (since it is defined as a line passing through two~$k$-points or a pair of conjugated points) and it intersects the toric boundary transversally. 
Hence, the motivic contribution equals $\qinv{1}$ (see Equation~\eqref{eq:mult}, or \cite{JPP23} for detailed computations). 
Thus, the motivic multiplicity of the tropical curve $\Gamma$ is trivially multiplicative as~$\mult(\Gamma)=\mult(\Gamma^{\text{t}})\cdot\mult(\Gamma_{\text{b}})$. 
Therefore, we have that the difference of the contributions of the tropical curves in this case vanishes since the tropical curves contributing to each count have the same bottom part. Thus,
\begin{align*}
&\sum_{\substack{\Gamma\supset \overline{\mathcal{P}}\\\deg(\Gamma^{\text{t}})=\Delta_1}}\mkern-15mu\mult(\Gamma)-\mkern-15mu\sum_{\substack{\Gamma\supset \overline{\mathcal{P}'}\\\deg(\Gamma^{\text{t}})=\Delta_1}}\mkern-15mu\mult(\Gamma)
\\&\quad=\left(\sum_{
\substack{\Gamma^{\text{t}}\supset\{\overline{p},\overline{p}'\}\\\deg(\Gamma^{\text{t}})=\Delta_1}}\mkern-15mu\mult(\Gamma^{\text{t}})
-\mkern-15mu\sum_{
\substack{\Gamma^{\text{t}}\supset\{\overline{q}\}\\\deg(\Gamma^{\text{t}})=\Delta_1}}\mkern-15mu\mult(\Gamma^{\text{t}})\right)\cdot
\mkern-15mu\sum_{
\substack{\Gamma_{\text{b}}\supset \overline{\mathcal{P}}\cap\overline{\mathcal{P}'} \\\deg(\Gamma_{\text{b}})=\overline{\Delta}}}\mkern-15mu\mult(\Gamma_{\text{b}})
\\&\quad=\left(\qinv{1}-\qinv{1}\right)\cdot
\mkern-15mu\sum_{
\substack{\Gamma_{\text{b}}\supset \overline{\mathcal{P}}\cap\overline{\mathcal{P}'} \\\deg(\Gamma_{\text{b}})=\overline{\Delta}}}\mkern-15mu\mult(\Gamma_{\text{b}})=0.
\end{align*}

\begin{figure}
\centering
\begin{tikzpicture}[line cap=round,line join=round,>=triangle 45,x=1.0cm,y=1.0cm,scale=0.8]
\clip(0,0) rectangle (5,4);
\draw [line width=2.pt] (0.,2.)-- (4.,2.);
\draw [line width=2.pt] (4.,0.)-- (4.,2.);
\draw [line width=2.pt] (4.,2.)-- (5.,3.);
\draw [line width=2.pt] (0.,3.)-- (3.,3.);
\draw [line width=2.pt] (3.,3.)-- (3.,0.);
\draw [line width=2.pt] (3.,3.)-- (4.,4.);
\draw [fill=red] (1.,3.) circle (2.5pt);
\draw(1.14,3.37) node {$\overline{p}$};
\draw [fill=RoyalBlue] (2.,2.) circle (2.5pt);
\draw (2.14,2.37) node {$\overline{p}'$};
\end{tikzpicture}
\hspace{0.5cm}
\begin{tikzpicture}[line cap=round,line join=round,>=triangle 45,x=1.0cm,y=1.0cm,scale=0.8]
\clip(0,0) rectangle (5,4);

\draw [line width=2.pt] (0.,2.5)-- (4.,2.5);
\draw [line width=2.pt] (4.,0.)-- (4.,2.5);
\draw [line width=2.pt] (4.,2.5)-- (5.,3.5);
\draw [line width=2.pt] (0.,2.5)-- (3.,2.5);
\draw [line width=2.pt] (3.,2.5)-- (3.,0.);
\draw [line width=2.pt] (3.,2.5)-- (4.5,4.);
\draw [fill=Purple] (1.5,2.5) circle (3pt);
\draw (1.5,2.1) node {$\overline{q}$};
\draw (0.3,2.75) node {\small$2$};
\end{tikzpicture}
    \caption{Tropical curves of degree $\Delta_2$ passing through two $k$-points or a double point defined over a quadratic extension $k(\sqrt{c})$, intersecting the bottom divisor in two $k$-points.}
    \label{fig:beta}
\end{figure}
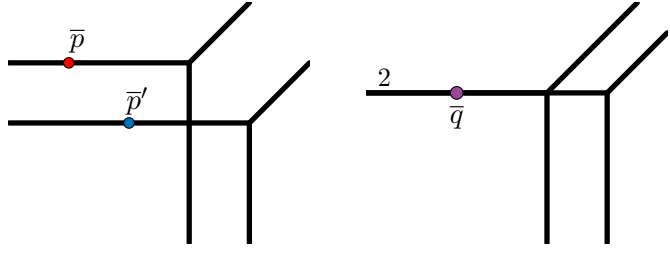

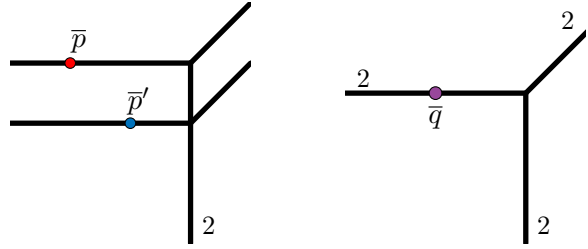
\begin{figure}
\centering
\begin{tikzpicture}[line cap=round,line join=round,>=triangle 45,x=1.0cm,y=1.0cm,scale=0.8]
\clip(0,0) rectangle (4,4);
\draw [line width=2.pt] (0.,2.)-- (3.,2.);
\draw [line width=2.pt] (3.,2.)-- (4.,3.);
\draw [line width=2.pt] (0.,3.)-- (3.,3.);
\draw [line width=2.pt] (3.,3.)-- (3.,0.);
\draw [line width=2.pt] (3.,3.)-- (4.,4.);
\draw [fill=red] (1.,3.) circle (2.5pt);
\draw(1.14,3.37) node {$\overline{p}$};
\draw [fill=RoyalBlue] (2.,2.) circle (2.5pt);
\draw (2.14,2.37) node {$\overline{p}'$};
\draw (3.3,0.3) node {\small$2$};
\end{tikzpicture}
\hspace{1cm}
\begin{tikzpicture}[line cap=round,line join=round,>=triangle 45,x=1.0cm,y=1.0cm,scale=0.8]
\clip(0,0) rectangle (4,4);

\draw [line width=2.pt] (0.,2.5)-- (3.,2.5);
\draw [line width=2.pt] (3.,2.5)-- (3.,0.);
\draw [line width=2.pt] (3.,2.5)-- (4.5,4.);
\draw [fill=Purple] (1.5,2.5) circle (3pt);
\draw (1.5,2.1) node {$\overline{q}$};
\draw (0.3,2.75) node {\small$2$};
\draw (3.7,3.7) node {\small$2$};
\draw (3.3,0.3) node {\small$2$};
\end{tikzpicture}
    \caption{Tropical curves of degree $\Delta_2$ passing through two $k$-points or a double point defined over a quadratic extension $k(\sqrt{c})$, intersecting the bottom divisor in a double point.}
    \label{fig:beta2}
\end{figure}

The second case consists of tropical curves that pass through the special points as in Figure~\ref{fig:beta} or as in Figure~\ref{fig:beta2}. 
These subcases are mutually exclusive. They can be distinguished by the field of definition of the point~${r}\in{\mathcal{P}}\cap{\mathcal{P}'}$, where $\overline{r}$ is the point on the elevator adjacent to the upper floor. In the former case, the point ${r}$ is $k$-rational, while in the latter one it is defined over a quadratic extension $k(\sqrt{c'})$. 
Like above, cutting the plane $\R^2$ at a horizontal line just above $\overline{r}$  splits each tropical curve $\Gamma$ into two pieces: the top part $\Gamma^{\text{t}}$ and the bottom part $\Gamma_{\text{b}}$. The curve $\Gamma^{\text{t}}$ has degree~$\Delta_2$, while the curve $\Gamma_{\text{b}}$  has as degree the Newton polygon $\widetilde{\Delta}\coloneqq\operatorname{NP}(\operatorname{Bl}_p(X))$ of the blow up of $X$ at a $k$-point.

For the subcase as in Figure~\ref{fig:beta},
 if the curve $\Gamma^{\text{t}}$ passes through~$\overline{p}$ and $\overline{p'}$, it is reducible and consists of two $k$-lines. Each reducible component has degree $\Delta_1$ and $\Gamma^{\text{t}}$ has degree the Minkowski sum $2\Delta_1=\Delta_2$. Its local motivic multiplicity is $\qinv{1}$ and it is trivially multiplicative. There are two such possibilities for the curve $\Gamma^{\text{t}}$ by exchanging $\overline{p}$ and $\overline{p}'$.
Else, if $\Gamma^{\text{t}}$ passes through~$\overline{q}$, its local motivic multiplicity equals $\Tra[k(\sqrt{c}))/k]{\qinv{1}}$ and it is multiplicative, by Lemma~\ref{lm:betamult}.
Therefore, in this case we have that
\begin{align*}
&\sum_{\substack{\Gamma\supset \overline{\mathcal{P}}\\\deg(\Gamma^{\text{t}})=2\Delta_1}}\mkern-15mu\mult(\Gamma)-\mkern-15mu\sum_{\substack{\Gamma\supset \overline{\mathcal{P}'}\\\deg(\Gamma^{\text{t}})=\Delta_2}}\mkern-15mu\mult(\Gamma)
\\&\quad=\left(\sum_{
\substack{\Gamma^{\text{t}}\supset\{\overline{p},\overline{p}'\}\\\deg(\Gamma^{\text{t}})=2\Delta_1}}\mkern-15mu\mult(\Gamma^{\text{t}})
-\mkern-15mu\sum_{
\substack{\Gamma^{\text{t}}\supset\{\overline{q}\}\\\deg(\Gamma^{\text{t}})=\Delta_2}}\mkern-15mu\mult(\Gamma^{\text{t}})\right)\cdot
\mkern-15mu\sum_{
\substack{\Gamma_{\text{b}}\supset \overline{\mathcal{P}}\cap\overline{\mathcal{P}'} \\\deg(\Gamma_{\text{b}})=\widetilde{\Delta}}}\mkern-15mu\mult(\Gamma_{\text{b}})
%\\&\quad=\left(2\qinv{1}-\Tra[k(\sqrt{c}))/k]{\qinv{1}}\right)\cdot
%\mkern-15mu\sum_{\substack{\Gamma_{\text{b}}\supset \overline{\mathcal{P}}\cap\overline{\mathcal{P}'} \\\deg(\Gamma_{\text{b}})=\widetilde{\Delta}}}\mkern-15mu\mult(\Gamma_{\text{b}}).
\\&\quad=\left(2\cdot\qinv{1}-\Tra[k(\sqrt{c}))/k]{\qinv{1}}\right)\cdot N^{\mathbb{A}^1}_{\operatorname{Bl}_p (X),D-2E,k}(\overline{\mathcal{P}}\cap\overline{\mathcal{P}'}).
\end{align*}

For the subcase as in Figure~\ref{fig:beta2}, the vertical edge passes through the point $\overline{r}$.
We consider an additional horizontal cut along the line that passes through~$\overline{r}$.
This cutting splits the tropical curve $\Gamma$ into two pieces of the same degrees as before: the top part $\Gamma^{\text{t}}_{-}$ and the bottom part $\Gamma_{\text{b}}^-$. The difference is the gluing of the algebraic curves, along the toric divisor $\operatorname{Tor}_k(\sigma)$, happens at the $k(\sqrt{c'})$-point $r$. Hence, there is a refinement of the point condition at the boundary given by the trace of $\qinv{1}$ from $k(\sqrt{c'})$ to $k$. By Lemma~\ref{lm:betamult}, the multiplicity of the tropical curve $\Gamma$ splits as 
\[\mult(\Gamma)=\mult(\Gamma^{\text{t}}_{-})
%\cdot\Tra[k(\sqrt{c'})/k]{\qinv{1}}
\cdot\mult(\Gamma_{\text{b}}^-),\]
with the difference that both pieces have point conditions at the boundary.
In order to compute the contribution from $\Gamma^{\text{t}}_-$, we can enlarge our field of definition to $k(\sqrt{c'})$, since the point conditions at the boundary are defined in this field extension, and trace down to $k$. Remark that in $k(\sqrt{c'})$, the points $r$ and its Galois conjugate are rational, so up to an isomorphism of the toric boundary, the tropical curves split as in Figure~\ref{fig:beta}. The trace forms of the respective multiplicities are
\begin{align*}
\Tra[k(\sqrt{c'})/k]{2\cdot\qinv{1}}&=2\cdot\qinv{1}\cdot\Tra[k(\sqrt{c'})/k]{\qinv{1}},\\
    \Tra[k(\sqrt{c'})/k]{\Tra[k(\sqrt{c},\sqrt{c'})/k(\sqrt{c'})]{\qinv{1}}}&=\Tra[k(\sqrt{c'})/k]{\qinv{1}}\cdot\Tra[k(\sqrt{c})/k]{\qinv{1}},
\end{align*}
by Lemma~\ref{lm:betamult}.
Hence, in the ring $\GW(k)$, we have that the difference of the multiplicities equals
\[\left(2\cdot\qinv{1}-\Tra[k(\sqrt{c})/k]{\qinv{1}}\right)\cdot
\Tra[k(\sqrt{c'})/k]{\qinv{1}}.
\]
To conclude, let us observe that in this subcase, a tropical curve $\Gamma_{\text{b}}\supset\overline{\mathcal{P}}\cap\overline{\mathcal{P}'}$ contributing to $N^{\mathbb{A}^1}_{\operatorname{Bl}_p (X),D-2E,k}(\overline{\mathcal{P}}\cap\overline{\mathcal{P}'})$ has multiplicity 
\[\mult(\Gamma_{\text{b}})=\Tra[k(\sqrt{c'})/k]{\qinv{1}}\cdot\mult(\Gamma_{\text{b}}^-),\]
since the point condition given by $\overline{r}$ has been moved to the boundary and its contribution is multiplicative due to Lemma~\ref{lm:betamult}. Therefore, we have that in this subcase
\begin{align*}
&\sum_{\substack{\Gamma\supset \overline{\mathcal{P}}\\\deg(\Gamma^{\text{t}})=\Delta_2}}\mkern-15mu\mult(\Gamma)-\mkern-15mu\sum_{\substack{\Gamma\supset \overline{\mathcal{P}'}\\\deg(\Gamma^{\text{t}})=\Delta_2}}\mkern-15mu\mult(\Gamma)
\\&\quad=
\left(2\cdot\qinv{1}-\Tra[k(\sqrt{c})/k]{\qinv{1}}\right)\cdot
\Tra[k(\sqrt{c'})/k]{\qinv{1}}
\cdot
\mkern-15mu\sum_{
\substack{\Gamma_{\text{b}}\supset \overline{\mathcal{P}}\cap\overline{\mathcal{P}'} \\\deg(\Gamma_{\text{b}})=\widetilde{\Delta}}}\mkern-15mu\mult(\Gamma_{\text{b}}^-)
%\\&\quad=\left(2\qinv{1}-\Tra[k(\sqrt{c}))/k]{\qinv{1}}\right)\cdot
%\mkern-15mu\sum_{\substack{\Gamma_{\text{b}}\supset \overline{\mathcal{P}}\cap\overline{\mathcal{P}'} \\\deg(\Gamma_{\text{b}})=\widetilde{\Delta}}}\mkern-15mu\mult(\Gamma_{\text{b}}).
\\&\quad=\left(2\cdot\qinv{1}-\Tra[k(\sqrt{c}))/k]{\qinv{1}}\right)\cdot N^{\mathbb{A}^1}_{\operatorname{Bl}_p (X),D-2E,k}(\overline{\mathcal{P}}\cap\overline{\mathcal{P}'}).
\end{align*}

The two main cases are exhausting for vertically stretched configurations. Since the cases are mutually exclusive for each tropical curve, we have that the difference of the counts $N^{\mathbb{A}^1}_{X,D,k}(\mathcal{P})$ and $N^{\mathbb{A}^1}_{X,D,k}(\mathcal{P}')$ is the sum of the computed differences.\qed

Let us remark that in this proof we studied the local behaviour of the tropical curve around its topmost part. 
Since the motivic counts are invariants, the same proof works for configurations of points having points defined over arbitrary field extensions as long as there is at least one point defined at most over a quadratic extension to separate them in the vertically stretched configuration. A similar proof can be constructed for higher degree extensions by studying tropical curves with higher degree points, but the combinatorics and local motivic contribution may be more challenging.

\section{Motivic plane invariants of degree 4 over quadratic field extensions}
\label{sec:apps}

In this section, we apply the wall crossing formula of Theorem~\ref{thm:WCF} in order to compute a complete list of the motivic plane invariants of degree~$4$ for configurations of points defined over quadratic field extensions. 
From now on, we denote by $4$ in the degree the class of $4L$ where $L$ is the strict transform of a generic line in $\mathbb{P}^2$.
We denote by $\mathcal{P}_s$ a generic configuration of points with exactly $s$ pairs of conjugated points, defined over the quadratic extensions $k(\sqrt{c_1}), k(\sqrt{c_2}),\dots, k(\sqrt{c_s})$, respectively. We define
\[\beta_i\coloneqq\Tra[k(\sqrt{c_i})/k]{\qinv{1}}.\]

In order to apply Theorem~\ref{thm:WCF} iteratively, we  compute the invariants for the  blow up $\Sigma_1$ of $\mathbb{P}^2$ at a~$k$-point, of degree $4-2E$; for the blow up $\Bl$ of $\Sigma_1$ at a $k$-point, of degree $4-2E-2E'$; and lastly, using that the motivic invariant of the blow up  of $\Sigma_1$ at two~$k$-points, of de\-gree~$4-2E-2E'-2E''$ is $\qinv{1}$.
We compute these invariants in backward order with respect to the number of blown up points. 

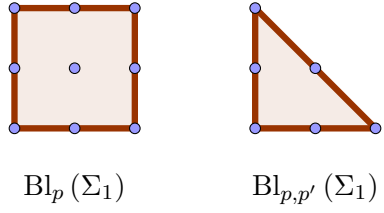
\begin{figure}
\centering
\definecolor{xdxdff}{rgb}{0.6,0.6,1.0}
\definecolor{ffffqq}{rgb}{1.0,1.0,0.0}
\definecolor{qqqqff}{rgb}{0.0,0.0,1.0}
\definecolor{zzttqq}{rgb}{0.6,0.2,0.0}

\begin{tikzpicture}[line cap=round,line join=round,>=triangle 45,x=1.0cm,y=1.0cm,scale=0.8]

\clip(11,-2) rectangle (19,5);
\fill[line width=2.4pt,color=zzttqq,fill=zzttqq,fill opacity=0.10000000149011612] (0.,0.) -- (4.,0.) -- (0.,4.) -- cycle;
\fill[line width=2.4pt,color=zzttqq,fill=zzttqq,fill opacity=0.10000000149011612] (6.,0.) -- (6.,2.) -- (8.,2.) -- (10.,0.) -- cycle;
\fill[line width=2.4pt,color=zzttqq,fill=zzttqq,fill opacity=0.10000000149011612] (12.,2.) -- (14.,2.) -- (14.,0.) -- (12.,0.) -- cycle;
\fill[line width=2.4pt,color=zzttqq,fill=zzttqq,fill opacity=0.10000000149011612] (16.,0.) -- (18.,0.) -- (16.,2.) -- cycle;
\draw [line width=2.4pt,color=zzttqq] (16.,0.)-- (18.,0.);
\draw [line width=2.4pt,color=zzttqq] (16.,0.)-- (16.,2.);
\draw [line width=2.4pt,color=zzttqq] (18.,0.)-- (16.,2.);

\draw [line width=2.4pt,color=zzttqq] (0.,0.)-- (4.,0.);
\draw [line width=2.4pt,color=zzttqq] (4.,0.)-- (0.,4.);
\draw [line width=2.4pt,color=zzttqq] (0.,4.)-- (0.,0.);
\draw [line width=2.4pt,color=zzttqq] (6.,0.)-- (6.,2.);
\draw [line width=2.4pt,color=zzttqq] (6.,2.)-- (8.,2.);
\draw [line width=2.4pt,color=zzttqq] (8.,2.)-- (10.,0.);
\draw [line width=2.4pt,color=zzttqq] (10.,0.)-- (6.,0.);
\draw [line width=2.4pt,color=zzttqq] (12.,2.)-- (14.,2.);
\draw [line width=2.4pt,color=zzttqq] (14.,2.)-- (14.,0.);
\draw [line width=2.4pt,color=zzttqq] (14.,0.)-- (12.,0.);
\draw [line width=2.4pt,color=zzttqq] (12.,0.)-- (12.,2.);

\draw [fill=xdxdff] (16.,0.) circle (2.5pt);
\draw [fill=xdxdff] (17.,0.) circle (2.5pt);
\draw [fill=xdxdff] (18.,0.) circle (2.5pt);
\draw [fill=xdxdff] (16.,1.) circle (2.5pt);
\draw [fill=xdxdff] (17.,1.) circle (2.5pt);
\draw [fill=xdxdff] (16.,2.) circle (2.5pt);

\draw [fill=xdxdff] (0.,0.) circle (2.5pt);
\draw [fill=xdxdff] (4.,0.) circle (2.5pt);
\draw [fill=xdxdff] (0.,4.) circle (2.5pt);
\draw [fill=xdxdff] (1.,0.) circle (2.5pt);
\draw [fill=xdxdff] (2.,0.) circle (2.5pt);
\draw [fill=xdxdff] (3.,0.) circle (2.5pt);
\draw [fill=xdxdff] (3.,1.) circle (2.5pt);
\draw [fill=xdxdff] (2.,1.) circle (2.5pt);
\draw [fill=xdxdff] (1.,1.) circle (2.5pt);
\draw [fill=xdxdff] (0.,1.) circle (2.5pt);
\draw [fill=xdxdff] (0.,2.) circle (2.5pt);
\draw [fill=xdxdff] (1.,2.) circle (2.5pt);
\draw [fill=xdxdff] (2.,2.) circle (2.5pt);
\draw [fill=xdxdff] (1.,3.) circle (2.5pt);
\draw [fill=xdxdff] (0.,3.) circle (2.5pt);
\draw [fill=xdxdff] (6.,0.) circle (2.5pt);
\draw [fill=xdxdff] (6.,2.) circle (2.5pt);
\draw [fill=xdxdff] (8.,2.) circle (2.5pt);
\draw [fill=xdxdff] (10.,0.) circle (2.5pt);
\draw [fill=xdxdff] (12.,2.) circle (2.5pt);
\draw [fill=xdxdff] (14.,2.) circle (2.5pt);
\draw [fill=xdxdff] (14.,0.) circle (2.5pt);
\draw [fill=xdxdff] (12.,0.) circle (2.5pt);
\draw [fill=xdxdff] (7.,2.) circle (2.5pt);
\draw [fill=xdxdff] (6.,1.) circle (2.5pt);
\draw [fill=xdxdff] (7.,1.) circle (2.5pt);
\draw [fill=xdxdff] (8.,1.) circle (2.5pt);
\draw [fill=xdxdff] (7.,0.) circle (2.5pt);
\draw [fill=xdxdff] (8.,0.) circle (2.5pt);
\draw [fill=xdxdff] (9.,0.) circle (2.5pt);
\draw [fill=xdxdff] (9.,1.) circle (2.5pt);
\draw [fill=xdxdff] (12.,1.) circle (2.5pt);
\draw [fill=xdxdff] (13.,1.) circle (2.5pt);
\draw [fill=xdxdff] (14.,1.) circle (2.5pt);
\draw [fill=xdxdff] (13.,2.) circle (2.5pt);
\draw [fill=xdxdff] (13.,0.) circle (2.5pt);

\draw (8.,-1.) node{$\Sigma_1$};
\draw (2.,-1.) node{$\mathbb{P}^2$};
\draw (13.,-1.) node{$\operatorname{Bl}_p\left(\Sigma_1\right)$};
\draw (17.,-1.) node{$\operatorname{Bl}_{p,p'}\left(\Sigma_1\right)$};
\end{tikzpicture}
    \caption{Lattice polygons corresponding to curves in $\operatorname{Bl}_p\left(\Sigma_1\right)$ of degree $4-2E-2E'$, and to curves in~$\operatorname{Bl}_{p,p'}\left(\Sigma_1\right)$ of degree $4-2E-2E'-2E''$, respectively.}
    \label{fig:polys2}
\end{figure}

Let us start by observing that the polygon associated to the blow up of~$\Sigma_1$ at two $k$-points of degree $4-2E-2E'-2E''$ is equal to $\Delta_2$, as in Figure~\ref{fig:polys2}. It follows from the fact that there is only one conic in $\mathbb{P}^2$ passing through $5$ points, either $k$-rationals or pairs of conjugated points over quadratic extensions, that the motivic invariant 
\[N^{\mathbb{A}^1}_{\operatorname{Bl}_{p,p'} (\Sigma_1),4-2E-2E'-2E'',k}=\qinv{1}.\]

Now, we use Theorem~\ref{thm:motcor} to compute the motivic invariant of $\Bl$ for a configuration that is $k$-rational. Applying the formula of Theorem~\ref{thm:WCF} to increase the number of points in the configuration defined over quadratic extensions yields
\begin{align*}
    N^{\mathbb{A}^1}_{\Bl,4-2E-2E',k}(\mathcal{P}_{s+1})=N^{\mathbb{A}^1}_{\Bl,4-2E-2E',k}(\mathcal{P}_s)+
%\Tra[k(\sqrt{c_{s+1}})/k]{\qinv{1}}
\beta_{s+1}
-2\cdot\qinv{1}.
\end{align*}
%\[N^{\mathbb{A}^1}_{\Bl,4-2E-2E',k}(\mathcal{P}_{s+1})=N^{\mathbb{A}^1}_{\Bl,4-2E-2E',k}(\mathcal{P}_s)+\Tra[k(\sqrt{c_{s+1}})/k]{\qinv{1}}-2\cdot\qinv{1}.\]
Therefore, the motivic invariants equal
\begin{align*}
    N^{\mathbb{A}^1}_{\Bl,4-2E-2E',k}(\mathcal{P}_0)&=2h+8\cdot\qinv{1},\\
N^{\mathbb{A}^1}_{\Bl,4-2E-2E',k}(\mathcal{P}_1)
&=2h+6\cdot\qinv{1}+\beta_1,\\
N^{\mathbb{A}^1}_{\Bl,4-2E-2E',k}(\mathcal{P}_2)
&=2h+4\cdot\qinv{1}+\beta_1+\beta_2,\\
N^{\mathbb{A}^1}_{\Bl,4-2E-2E',k}(\mathcal{P}_3)
&=2h+2\cdot\qinv{1}+\beta_1+\beta_2+\beta_3.\\
\end{align*}

The next step is to compute the motivic invariants of $\Bl$, using  Theorem~\ref{thm:motcor} for a configuration that is $k$-rational and applying the formula of Theorem~\ref{thm:WCF} as before. This formula yields
\begin{align*} N^{\mathbb{A}^1}_{\Sigma_1,4-2E,k}(\mathcal{P}_{s+1})=&N^{\mathbb{A}^1}_{\Sigma_1,4-2E,k}(\mathcal{P}_s)\\&+
\left(\beta_{s+1}-2\cdot\qinv{1}\right)\cdot N^{\mathbb{A}^1}_{\Bl,4-2E-2E',k}(\mathcal{P}_{s}).
\end{align*}
Therefore, the motivic invariants equal
\begin{align*}
N^{\mathbb{A}^1}_{\Sigma_1,4-2E,k}(\mathcal{P}_0)&=24h+48\cdot\qinv{1}, \\
N^{\mathbb{A}^1}_{\Sigma_1,4-2E,k}(\mathcal{P}_1)&=24h+32\cdot\qinv{1} +8\cdot\beta_1,\\
N^{\mathbb{A}^1}_{\Sigma_1,4-2E,k}(\mathcal{P}_2)&=24h+20\cdot\qinv{1} +6\cdot\beta_1+6\cdot\beta_2+\beta_1\beta_2,\\
N^{\mathbb{A}^1}_{\Sigma_1,4-2E,k}(\mathcal{P}_3)&=24h+12\cdot\qinv{1} +4\cdot\sum_{i=1}^3\beta_i+\sum_{1\leq i<j\leq 3}\beta_i\beta_j,\\
N^{\mathbb{A}^1}_{\Sigma_1,4-2E,k}(\mathcal{P}_4)&=24h+8\cdot\qinv{1} +2\cdot\sum_{i=1}^4\beta_i+\sum_{1\leq i<j\leq 4}\beta_i\beta_j,\\
\end{align*}
The final step is to compute the motivic invariants of $\mathbb{P}^2$, using  Theorem~\ref{thm:motcor} for an initial step and iterating the formula of Theorem~\ref{thm:WCF}
\[
N^{\mathbb{A}^1}_{\mathbb{P}^2,4,k}(\mathcal{P}_{s+1})=N^{\mathbb{A}^1}_{\mathbb{P}^2,4,k}(\mathcal{P}_s)+
\left(\beta_{s+1}-2\cdot\qinv{1}\right)\cdot N^{\mathbb{A}^1}_{\Sigma_1,4-2E,k}(\mathcal{P}_{s}).
\]
Therefore, the motivic invariants equal
\begin{align*}
N^{\mathbb{A}^1}_{\mathbb{P}^2,4,k}(\mathcal{P}_0)&=190h+240\cdot\qinv{1}, \\
N^{\mathbb{A}^1}_{\mathbb{P}^2,4,k}(\mathcal{P}_1)&=190h+144\cdot \qinv{1}+48\cdot \beta_1, \\
N^{\mathbb{A}^1}_{\mathbb{P}^2,4,k}(\mathcal{P}_2)&=190h+80\cdot \qinv{1}+32\cdot\left(\beta_1+\beta_2\right)+8\cdot\beta_1\beta_2, \\
N^{\mathbb{A}^1}_{\mathbb{P}^2,4,k}(\mathcal{P}_3)&=190h+40\cdot \qinv{1}+20\cdot\sum_{i=1}^3\beta_i+6\cdot\mkern-15mu\sum_{1\leq i<j\leq 3}\beta_i\beta_j+\beta_1\beta_2\beta_3, \\
N^{\mathbb{A}^1}_{\mathbb{P}^2,4,k}(\mathcal{P}_4)&=190h+16\cdot \qinv{1}+12\cdot\sum_{i=1}^4\beta_i+4\cdot
\mkern-15mu\sum_{1\leq i<j\leq 4}\beta_i\beta_j+
\mkern-15mu\sum_{1\leq i<j<l\leq 4}\beta_i\beta_j\beta_l, \\
N^{\mathbb{A}^1}_{\mathbb{P}^2,4,k}(\mathcal{P}_5)&=190h+8\cdot\sum_{i=1}^5\beta_i+2\cdot
\mkern-15mu\sum_{1\leq i<j\leq 5}\beta_i\beta_j+
\mkern-15mu\sum_{1\leq i<j<l\leq 5}\beta_i\beta_j\beta_l. \\
\end{align*}

Let us remark that the last line of each set of motivic invariants presented in this section can be seen as a polynomial in $\GW(k)\left[\bar{\beta}\right]$ that specialize to configurations with any number of $k$-points by setting any $\beta_i$ as the trace form $\Tra[k\times k/k]{\qinv{1}}=2\cdot \qinv{1}$, or equivalently, choosing $d_i$ as a square in~$k^*$. We opted for the current  presentation to allow for a parallel with other refined invariants, like the Block-Göttsche invariants introduced in~\cite{BG14}, whose polynomial properties for invariants counting rational curves are studied in~\cite{BJP}.

The relation given by Theorem~\ref{thm:WCF} implies that $N^{\mathbb{A}^1}_{X,D,k}(\mathcal{P})$ is a monic polynomial of degree given by the number of interior points $g$ of the polygon associated to $(X,D)$, since after $g$ blow ups we have an invariant determined by its real and complex realization by Theorem~\ref{thm:motcor}.

\bibliographystyle{plain}
\bibliography{bibliographie.bib}

\end{document}